\documentclass[12pt]{article}
\usepackage{amssymb,amscd}
\usepackage{graphicx}

\newtheorem{conj}{Conjecture}[section]
\newtheorem{theo}[conj]{Theorem}
\newtheorem{rem}[conj]{Remark}

\newtheorem{defin}[conj]{Definition}
\newtheorem{prop}[conj]{Proposition}
\newtheorem{cor}[conj]{Corollary}
\newtheorem{lema}[conj]{Lemma}
\newtheorem{problem}[conj]{Problem}
\newtheorem{exam}[conj]{Example}

\begin{document}
\date{October 18, 2015}
\title{\Large Multiple chessboard complexes\\ and the colored Tverberg problem}

\author{{\large Du\v sko Joji\'c }
\\ {Faculty of Science, University of Banja Luka}\\[2mm]
  {\large Sini\v sa T. Vre\' cica\thanks{Supported by the Ministry
for Science and Technology of Serbia, Grant 174034.}}
\\ {Faculty of Mathematics, University of Belgrade}\\[2mm]
  {\large  Rade T. \v Zivaljevi\' c${}^\ast$}
 \\ {Mathematics Institute
 SASA, Belgrade}}
\maketitle 

\begin{abstract}
Following D.B. Karaguezian, V. Reiner, and M.L. Wachs (Matching
Complexes, Bounded Degree Graph Complexes, and Weight Spaces of
$GL$-Complexes, Journal of Algebra 2001) we study the connectivity
degree and shellability of multiple chessboard complexes. Our
central new results (Theorems~\ref{thm:main} and
\ref{thm:shelling-main}) provide sharp connectivity bounds
relevant to applications in Tverberg type problems where multiple
points of the same color are permitted. These results also provide
a foundational work for the new results of Tverberg-van
Kampen-Flores type, as announced in the forthcoming paper
\cite{jvz}.

\end{abstract}

\section{An overview and motivation}

Chessboard complexes and their generalizations belong to the class
of most studied graph complexes, with numerous applications in and
outside combinatorics \cite{A04, blvz, bmz, FH, Garst, Jo-book,
krw, M, ShWa, VZ94, zv92, Ziegler, Z04}.

\medskip
The {\em connectivity degree} of a simplicial complex was selected
in \cite[Chapter 10]{Jo-book} as one of the five most important
and useful parameters in the study of simplicial complexes of
graphs. Following \cite{krw} we study the connectivity degree of
{\em multiple chessboard complexes} (Section~\ref{sec:multiple})
and their generalizations (Section~\ref{sec:generalized}). Our
first central result is Theorem~\ref{prva} which improves the
$2$-dimensional case of \cite[Corollary~5.2.]{krw} and reduces to
the $2$-dimensional case of \cite[Theorem~3.1.]{blvz} in the case
of standard chessboard complexes.

\medskip
Perhaps it is worth emphasizing that our methods allow us to
obtain sharp bounds relevant to applications in Tverberg and van
Kampen-Flores type problems (see Section~\ref{thm:application} and
\cite{jvz}). Moreover, the focus in \cite{krw} is on the homology
with the coefficients in a field and multidimensional chessboard
complexes while our results are homotopical and apply to
$2$-dimensional chessboard complexes.

\medskip
High connectivity degree is sometimes a consequence of the
shellability of the complex (or one of its skeletons), see
\cite{zie} for an early example in the context of chessboard
complexes. Theorem~\ref{thm:shelling-main} provides a sufficient
condition which guarantees the shellability of multiple chessboard
complexes and yields another proof of Theorem~\ref{thm:main}. The
construction of the shelling offers a novel point of view on this
problem and seems to be new and interesting already in the case of
standard chessboard complexes.

\medskip
Among the initial applications of the new connectivity bounds
established by Theorem~\ref{prva} is a result of colored Tverberg
type where multiple points of the same color are permitted
(Theorem~\ref{thm:application} in Section~\ref{sec:application}).
After the first version of our paper was submitted to the arXiv we
were kindly informed by G\"{u}nter Ziegler that
Theorem~\ref{thm:application} is implicit in their recent work
(see \cite{BFZ}, Theorem~4.4 and the remark following the proof of
Lemma~4.2.).

\medskip
Other, possibly more far reaching applications of
Theorems~\ref{thm:main} and \ref{thm:shelling-main} to theorems of
Tverberg-van Kampen-Flores type are announced in \cite{jvz}. This
provides new evidence that the chessboard complexes and their
generalizations are a natural framework for constructing
configuration spaces relevant to Tverberg type problems and
related problems about finite sets of points in Euclidean spaces.

\medskip\noindent {\small
{\bf Caveat:} Most if not all simplicial complexes in this paper
are visualized in a rectangular chessboard $[m]\times [n]$. The
reader is free to choose either the Cartesian or the matrix
enumeration of squares (where $(1,1)$ is the lower left corner in
the first and the upper left corner in the second). This should
not generally affect the reading of the paper and a little care is
needed only when interpreting the Figure~\ref{fig:sah-1},
revealing our slight inclination towards the Cartesian notation. }

\subsection{Colored Tverberg problems}
\label{sec:intro-1}

`Tverberg problems' is a common name for a class of theorems and
conjectures about finite sets of points (point clouds) in
$\mathbb{R}^{d}$. We start with a brief introduction into this
area of topological combinatorics emphasizing, in the spirit of
\cite{TV-bundle} and \cite{vz11}, a graphical or diagrammatic
((\ref{eqn:summary-Tverberg})-(\ref{eqn:(4)})) presentation of
these results. The reader is referred to \cite{Ziegler},
\cite[Section~14.4]{vz11}, \cite{Z04}, and \cite{M} for more
complete expositions of these problems and the history of the
whole area.

\medskip
The Tverberg theorem \cite{Tve} claims that every set $K\subset
\mathbb{R}^d$ with $(d+1)(q-1)+1$ elements can be partitioned $K =
K_1\cup\ldots\cup K_q$ into $q$ nonempty, pairwise disjoint
subsets $K_1, \ldots, K_q$ such that the corresponding convex
hulls have a nonempty intersection:
\begin{equation}
\bigcap_{i=1}^q {\rm conv}(K_i)\neq\emptyset .
\end{equation}
Following \cite{bss} it can be reformulated as the statement that
for each linear (affine) map $f : \Delta^D
\stackrel{a}{\longrightarrow} \mathbb{R}^d$ ($D =(d+1)(q-1)$)
there exist $q$ nonempty disjoint faces $\Delta_1,\ldots,
\Delta_q$ such that $f(\Delta_1)\cap\ldots\cap
f(\Delta_q)\neq\emptyset$. This form of Tverberg's result can be
summarized as follows,
\begin{equation}\label{eqn:summary-Tverberg}
(\Delta^{(d+1)(q-1)} \stackrel{a}{{\longrightarrow}} \mathbb{R}^d)
\Rightarrow (q- {\rm intersection} ).
\end{equation}
Here we tacitly assume that the faces intersecting in the image
are always vertex disjoint. The letter ``a'' over the arrow means
that the map is affine and its absence indicates that it can be an
arbitrary continuous map.

\medskip
The following four statements are illustrative for results of
`colored Tverberg type'.
\begin{equation}\label{eqn:(1)}
(K_{3,3} \longrightarrow {\mathbb R}^2) \Rightarrow (2- {\rm
intersection})
\end{equation}
\vspace{-7mm}
\begin{equation}\label{eqn:(2)}
(K_{3,3,3} \stackrel{a}{{\longrightarrow}} {\mathbb R}^2)
\Rightarrow (3- {\rm intersection})
\end{equation}
\begin{equation}\label{eqn:(3)}
(K_{5,5,5} \longrightarrow {\mathbb R}^3) \Rightarrow (3- {\rm
intersection})
\end{equation}
\begin{equation}\label{eqn:(4)}
(K_{4,4,4,4} \longrightarrow {\mathbb R}^3) \Rightarrow (4- {\rm
intersection})
\end{equation}

\medskip\noindent $K_{t_1,t_2,\ldots,
t_k}=[t_1]\ast[t_2]\ast\ldots\ast[t_k]$ is by definition the
complete multipartite simplicial complex obtained as a join of
$0$-dimensional complexes (finite sets). By definition the
vertices of this complex are naturally partitioned into groups of
the same `color'. For example $K_{p,q}=[p]\ast[q]$ is the complete
bipartite graph obtained by connecting each of $p$ `red vertices'
with each of $q$ `blue vertices'. The simplices of
$K_{t_1,t_2,\ldots, t_k}$ are called multicolored sets or {\em
rainbow simplices} and its dimension is $k_{-}=k-1$. (We
systematically use the abbreviation $a_{-}:= a-1$ to emphasize
that $a_{-}$ is the dimension of a non-degenerate simplex with $a$
vertices.)

\medskip
The implication (\ref{eqn:(1)}) says that for each continuous map
$\phi : K_{3,3}\rightarrow \mathbb{R}^2$ there always exist two
vertex-disjoint edges which intersect in the image. In light of
the Hanani-Tutte theorem this statement is equivalent to the
non-planarity of the complete bipartite graph $K_{3,3}$. The
implication (\ref{eqn:(2)}) is an instance of a result of
B\'{a}r\'{a}ny and Larman \cite{BL}. It says that each collection
of nine points in the plane, evenly colored by three colors, can
be partitioned into three multicolored or `rainbow triangles'
which have a common point.

\medskip
Note that a $9$-element set $C\subset \mathbb{R}^2$ which is
evenly colored by three colors, can be also described by a map
$\alpha : [3]\sqcup [3] \sqcup [3] \rightarrow \mathbb{R}^2$ from
a disjoint sum of three copies of $[3]$. In the same spirit an
affine map $\phi : K_{3,3,3}\stackrel{a}{\longrightarrow}
\mathbb{R}^2$ parameterizes not only the colored set itself but
takes into account from the beginning that some simplices
(multicolored or rainbow simplices) play a special role.

\medskip
A similar conclusion has statement (\ref{eqn:(3)}) which is a
formal analogue of the statement (\ref{eqn:(1)}) in dimension $3$.
It is an instance of a result of Vre\' cica and \v{Z}ivaljevi\' c
\cite{VZ94}, which claims the existence of three intersecting,
vertex disjoint rainbow triangles in each constellation of $5$
red, $5$ blue, and $5$ green stars in the $3$-space. A consequence
of this result is that $K_{5,5,5}$ is strongly non-embeddable in
$\mathbb{R}^3$ in the sense that there always exists a triple
point in the image.

\medskip
Finally (\ref{eqn:(4)}) is an instance of the celebrated result of
Blagojevi\'{c}, Matschke, and Ziegler \cite[Corollary~2.4]{bmz}
saying that $4$ intersecting, vertex disjoint rainbow tetrahedra
in $\mathbb{R}^3$ will always appear if we are given sixteen
points, evenly colored by four colors.

\begin{rem}{\rm Both statements (\ref{eqn:(4)}) and
(\ref{eqn:(3)}) are instances of results of colored Tverberg type.
There is an important difference between them however, and this is
the reason why they are referred to as Type~A and Type~B colored
Tverberg theorems in the Handbook of discrete and computational
geometry \cite[Chapter~14]{Z04}. Both results are optimal in the
sense that in the cases where they apply they provide the best
bounds possible.

}
\end{rem}

\subsection{General colored Tverberg theorems}
\label{sec:general-color}

From the point of view of results exhibited in
Section~\ref{sec:intro-1} it is quite natural to ask the following
general question.

\begin{problem}\label{prob:general-1}
For given integers $r, k, d$ determine the smallest $t = T(r,k,d)$
such that,
\begin{equation}\label{eqn:(gen-1)}
(K_{t,t,\ldots,t} \longrightarrow {\mathbb R}^d) \Rightarrow
(r-{\rm intersection})
\end{equation}
where $K_{t,t,\ldots,t} = K_{t,t,\ldots,t; k} = [t]^{\ast (k+1)}$
is the join of $k+1$ copies of $[t]$.
\end{problem}
The latest developments \cite{Ziegler, bmz} showed the importance
of the following even more general, non-homogeneous version of
Problem~\ref{prob:general-1}

\begin{problem}\label{prob:general-2}
For given integers $r, k, d$ determine when a sequence
$\mathfrak{t}=(t_0,t_1,\ldots, t_k)$ yields the implication
\begin{equation}\label{eqn:(gen-2)}
(K_{t_0,t_1,\ldots,t_k} \longrightarrow {\mathbb R}^d) \Rightarrow
(r-{\rm intersection})
\end{equation}
where $K_{\mathfrak{t}}=K_{t_0,t_1,\ldots,t_k} = [t_0]\ast
[t_1]\ast \ldots \ast [t_k]$.
\end{problem}

Historically the first appearance of the colored Tverberg problem
is the question of B\' ar\' any and Larman \cite{BL}. It is
related to the case $k = d$ of Problem~\ref{prob:general-1} i.e.\
to the case when the dimension $k$ of top-dimensional rainbow
simplices is equal to the dimension $d$ of the ambient Euclidean
space. This case is referred to in \cite{Z04} as the Type A of the
colored Tverberg problem. The Type B of colored Tverberg problem,
corresponding to the case $k < d$ of Problem~\ref{prob:general-1}
is introduced in \cite{VZ94}, see also \cite{User1,User2,Z04}.

\medskip
The following two theorems are currently the most general known
results about the invariants $T(r,k,d)$. The first is the recent
Type A statement due to Blagojevi\'{c}, Matschke, and Ziegler
\cite{bmz} who improved the original Type A colored Tverberg
theorem of Vre\' cica and \v Zivaljevi\' c \cite{zv92}. The second
is a Type B statement proved by Vre\' cica and \v Zivaljevi\' c in
\cite{VZ94}. Both results are exact in the sense that in the cases
where they apply they both evaluate the exact value of the
function $T(r,k,d)$.

\begin{theo}{\rm (\cite{bmz})}\label{thm:type-A}
If  $r+1$ is a prime number then $T(r,d,d)= r$.
\end{theo}
\begin{theo}{\rm (\cite{VZ94})}\label{thm:type-B}
If  $k\le d$ and $r$ is a prime such that $2\leq r \leq d/(d-k)$
then
$$T(r,k,d)= 2r-1.$$
\end{theo}
The condition $r \leq d/(d-k)$ in the second statement cannot be
removed if $r\geq 2$. It does not exist in the Type A statement
(Theorem~\ref{thm:type-A}) so it is a characteristic,
distinguishing feature of the Type B colored Tverberg theorem
(Theorem~\ref{thm:type-B}).

\subsection{Chessboard complexes and colored Tverberg problem}
\label{sec:general-scheme}

Here we briefly outline, following the original sources
\cite{zv92,VZ94} and a more recent exposition given in
\cite{vz11}, how the so called {\em chessboard complexes}
naturally arise in the context of the colored Tverberg problem.
For notation and a more systematic exposition of these and related
facts the reader is referred to \cite{M,User1,User2,Z04}.

\medskip
Given a map $f : K \rightarrow \mathbb{R}^d$ (as in examples from
Sections~\ref{sec:intro-1} and \ref{sec:general-color}) we want to
find $r$ nonempty, vertex disjoint faces $\sigma_1,\ldots,
\sigma_r$ of $K$ such that $f(\sigma_1)\cap\ldots\cap
f(\sigma_r)\neq\emptyset$. For this reason we consider the induced
map $F : K^{\ast r}_{\Delta}\rightarrow (\mathbb{R}^d)^{\ast r}$
from the deleted join (see \cite[Sections~5.5 and 6.3]{M}) of $r$
copies of $K$ to the $r$-fold join of $\mathbb{R}^d$ and observe
that it is sufficient to show that ${\rm Image}(F)\cap D \neq
\emptyset$ where $D \subset (\mathbb{R}^d)^{\ast r}$ is the
diagonal subspace of the join. Assuming the contrary we obtain a
$S_r$-equivariant map $F' : K^{\ast r}_{\Delta}\rightarrow
(\mathbb{R}^d)^{\ast r}_\Delta$ where $(\mathbb{R}^d)^{\ast
r}_\Delta$ is the $r$-fold deleted join of $\mathbb{R}^d$
\cite[Section~6.3]{M}. It is not difficult to show that
$(\mathbb{R}^d)^{\ast r}_\Delta$ has the $S_r$-homotopy type of
the unit sphere $S(W_r^{\oplus (d+1)})$ where $W_r$ is the
$(r-1)$-dimensional, standard (real) representation of $S_r$.

\medskip
If $K = K_{t_1,\ldots, t_k} = [t_1]\ast\ldots\ast [t_k]$ (as in
examples from Sections~\ref{sec:intro-1} and
\ref{sec:general-color}) then,
\begin{equation}
K^{\ast r}_{\Delta} \cong ([t_1]\ast\ldots\ast [t_k])^{\ast
r}_\Delta \cong [t_1]^{\ast r}_\Delta\ast\ldots \ast[t_k]^{\ast
r}_\Delta = \Delta_{r,t_1}\ast\ldots\ast\Delta_{r,t_k}
\end{equation}
where $\Delta_{m,n}$ is the so called {\em chessboard complex},
defined as the simplicial complex of all non-taking rook
placements on a $m\times n$ `chessboard'.

\medskip
The upshot of this sequence of reductions is that the implication
(\ref{eqn:(gen-2)}) (Problem~\ref{prob:general-2}) is a
consequence of a Borsuk-Ulam type result claiming that here does
not exits a $S_r$-equivariant map
\begin{equation}
F : \Delta_{r,t_1}\ast\ldots\ast\Delta_{r,t_k} \longrightarrow
S(W_r^{\oplus (d+1)}).
\end{equation}
In particular Theorems~\ref{thm:type-A} and \ref{thm:type-B} are
both reduced to the question of non-existence of $S_r$-equivariant
maps of spaces where the source space is a join of chessboard
complexes,
\begin{equation}
F : (\Delta_{r,r-1})^{\ast d}\ast \Delta_{r,1} \rightarrow
S(W_r^{\oplus d}).
\end{equation}
\begin{equation}
F : (\Delta_{r,2r-1})^{\ast k} \rightarrow S(W_r^{\oplus d}).
\end{equation}

\subsection{Multiple chessboard complexes}
\label{sec:multiple}

It is quite natural to apply the scheme outlined in
Section~\ref{sec:general-scheme} to some other simplicial
complexes aside from $K_{t_1,\ldots, t_k}$.

Let $[t]^{(p)}=\{A\subset [t] \mid \vert A\vert\leq p\}$ be the
collection of all subsets of $[t]=\{1,\ldots, t\}$ of size at most
$p$. As a simplicial complex $[t]^{(p)}$ is the $(p-1)$-skeleton
of the simplex spanned by $[t]$. Let
$K^{\mathfrak{p}}_{{\mathfrak{t}}}=K^{p_1,\ldots, p_k}_{t_1,
\ldots,t_k} = [t_1]^{(p_1)}\ast\ldots\ast [t_k]^{(p_k)}$.

\medskip
By applying the same procedure as in
Section~\ref{sec:general-scheme} to the complex
$K^{\mathfrak{p}}_{\mathfrak{t}}=K^{p_1,\ldots, p_k}_{t_1,
\ldots,t_k}$ we have,
\begin{equation}
(K^{\mathfrak{p}}_{{\mathfrak{t}}})^{\ast r}_{\Delta} \cong
([t_1]^{(p_1)}\ast\ldots\ast [t_k]^{(p_k)})^{\ast r}_\Delta \cong
([t_1]^{(p_1)})^{\ast r}_\Delta\ast\ldots
\ast([t_k]^{(p_k)})^{\ast r}_\Delta =
\Delta_{r,t_1}^{1,p_1}\ast\ldots\ast\Delta_{r,t_k}^{1,p_k}
\end{equation}
where the generalized chessboard complexes $\Delta_{r,t}^{1,p} =
[t]^{(p)})^{\ast r}_\Delta$ make their first appearance in this
paper. Formally they are introduced in greater generality in
Section~\ref{sec:generalized}.

\subsection{Bier spheres as multiple chessboard complexes}

One of the novelties in the proof of Theorem~\ref{thm:type-A}
\cite{bmz}, which is particularly visible in the `mapping degree'
proof \cite{vz11} and \cite{BMZ-2}, is the use of the
(pseudo)manifold structure of the chessboard complex
$\Delta_{r,r-1}$. Here we observe that an important subclass of
combinatorial spheres (Bier spheres) arise as multiple chessboard
complexes. As shown in Example~\ref{exam:Bier} all Bier spheres
can be incorporated into this scheme if we allow even more general
chessboard complexes.

\medskip
Recall \cite[Chapter~5]{M} that the Bier sphere $Bier_m(K)$,
associated to a simplicial complex $K\subset 2^{[m]}$, is the
`disjoint join' $K\ast K^\circ$ of $K$ and its combinatorial
Alexander dual $K^\circ$. The reader is referred to
Definitions~\ref{def:general-chessboard}-\ref{def:special-2} in
Section~\ref{sec:generalized} for the definition of generalized
chessboard complexes $\Delta_{m,2}^{\mathbf{k}; \mathbf{l}},
\Delta_{m,2}^{m-2,1(1)}$ and their relatives.

\begin{prop}\label{prop:Bier-1}
Suppose that $K = [m]^{\leq p}$ is the simplicial complex of all
subsets of $[m]$ of size at most $p$ and let ${\rm Bier}_m(K)$ be
the associated Bier sphere. Then
\[
S^{m-2}\cong {\rm Bier}_m(K) = {\rm Bier}_m([m]^{\leq p}) \cong
\Delta_{m,2}^{\mathbf{k}; \mathbf{l}}
\]
where $\mathbf{k} = (p, m-p-1)$ and $l_1=\ldots=l_m=1$, in
particular
\[
\Delta_{2p+1,2}^{p,1}\cong {\rm Bier}_{2p+1}([2p+1]^{\leq p})\cong
S^{2p-1} \quad \mbox{\rm and} \quad \Delta_{m,2}^{m-2,1(1)}\cong
{\rm Bier}_m([m])\cong S^{m-2}.
\]
\end{prop}

\section{Generalized chessboard complexes}
\label{sec:generalized}

The classical chessboard complex $\Delta_{m,n}$  \cite{blvz} is
often visualized as the simplicial complex of non-taking rook
placements on a $(m\times n)$-chessboard. In particular its
vertices $Vert(\Delta_{m,n}) = [m]\times [n]$ are elementary
squares in a chessboard which has $n$ rows of size $m$ (here we
use Cartesian rather than matrix presentation of the chessboard).

\medskip
The complex $\Delta_{m,n}$  can be also described as the {\em
matching complex} of the complete bipartite graph $K_{m,n}$. In
this incarnation its vertices correspond to all edges of the graph
$K_{m,n}$ and a collection of edges determine a simplex if and
only if it is a matching in $K_{m,n}$, see \cite{blvz} or
\cite{Jo-book}. As we have already seen in
Section~\ref{sec:general-scheme} the complex $\Delta_{m,n}$  can
be also described as the $n$-fold $2$-deleted join of the
$0$-dimensional skeleton of the $(m-1)$-dimensional simplex
$$((\sigma^{m-1})^{(0)})_{\Delta (2)}^{*n}.$$
Here, the $n$-fold $(q+1)$-deleted join of the complex $K$,
denoted by $K_{\Delta (q+1)}^{*n}$, is a subcomplex of $K^{*n}$,
the $n$-fold join of the complex $K$, consisting of joins of
$n$-tuples of simplices from $K$ such that the intersection of any
$q+1$ of them is empty. (In particular the $2$-deleted join
$K_{\Delta (2)}^{*n} = K_{\Delta}^{*n}$ is the usual deleted join
of $K$.)

The $0$-dimensional skeleton of the $(m-1)$-dimensional simplex
and the $m$-fold $2$-deleted join of a point are both identified
as the sets of $m$ points. It follows that,
$$
\Delta_{m,n} = ((\sigma^{m-1})^{(0)})_{\Delta (2)}^{*n} \cong
([\ast ]_{\Delta (2)}^{*m})_{\Delta (2)}^{*n}.
$$
This is precisely the description of $\Delta_{m,n}$ that appeared
in the original approach to the Colored Tverberg theorem in
\cite{zv92}.

\medskip
In this paper we allow multicolored simplices to have more (say
$p$) vertices of the same color, so we consider a {\em generalized
chessboard complex} which is the $n$-fold $(q+1)$-deleted join of
the $(p-1)$-dimensional skeleton of the $(m-1)$-dimensional
simplex, i.e. the complex
\begin{equation}\label{eqn:general}
\Delta_{m,n}^{p,q} := ((\sigma^{m-1})^{(p-1)})_{\Delta
(q+1)}^{*n}=([\ast ]_{\Delta (p+1)}^{*m})_{\Delta (q+1)}^{*n}
\end{equation}

As before, the vertices of this simplicial complex correspond to
the squares on the $m\times n$ chessboard and simplices correspond
to the collections of vertices so that at most $p$ of them are in
the same row, and at most $q$ of them are in the same column. As
indicated in (\ref{eqn:general}) we denote this simplicial complex
by $\Delta_{m,n}^{p,q}$  so in particular $\Delta_{m,n}=
\Delta_{m,n}^{1,1}$.

\begin{rem}\label{rem:caveat}{\rm The meaning of parameters
($m, p; n, q$)) in the complex $\Delta_{m,n}^{p,q}$ can be
memorized as follows. The parameters $m$ and $p$ both apply to the
rows ($m$ as the row-length of the chessboard $[m]\times [n]$ and
$p$ as the maximum number of rooks allowed in each row).
Similarly, the parameters $n$ and $q$ are associated to columns
($[n]$ is the column-height of $[m]\times [n]$ while $q$
prescribes the largest number of rooks in each of the columns). A
similar interpretation can be given in the case of more general
chessboard complexes ((\ref{eqn:2-collections}) and
(\ref{eqn:defin-hv})). }
\end{rem}

\begin{rem}{\rm The higher dimensional
analogues of complexes $\Delta_{m,n}^{p,q}$ were introduced and
studied in \cite{krw} and our particular interest in the
generalized Tverberg-type problems is the reason why in this paper
we focus on the two dimensional case. }
\end{rem}

\subsection{Complexes $\Delta_{m,n}^{\mathcal{K}, \mathcal{L}}$}

Both for heuristic and technical reasons we consider even more
general chessboard complexes based on the $(m\times
n)$-chessboard.  The following definition provides an ecological
niche (and a summary of notation) for all these complexes.

\begin{defin}\label{def:general-chessboard}
Let $\mathcal{K} = \{K_i\}_{i=1}^n$ and $\mathcal{L} =
\{L_j\}_{j=1}^m$ be two labelled collections of simplicial
complexes where $Vert(K_i) = [m]$ for each $i\in [n]$ and
$Vert(L_j) = [n]$ for each $j\in [m]$. Define,
\begin{equation}\label{eqn:2-collections}
\Delta_{m,n}^{\mathcal{K}, \mathcal{L}} =
\Delta_{m,n}(\mathcal{K}, \mathcal{L})
\end{equation}
as the complex of all subsets (rook-placements) $A\subset
[m]\times [n]$ such that $\{i\in [m]\mid (i,j)\in A \}\in K_j$ for
each $j\in [n]$ and $\{j\in [n]\mid (i,j)\in A \}\in L_i$ for each
$i\in [m]$.
\end{defin}
\begin{exam}\label{exam:Bier}{\rm
Generalizing Proposition~\ref{prop:Bier-1} we observe that the
general Bier sphere $Bier_m(K)$ arises as the complex
$\Delta_{m,2}^{\mathcal{K}, \mathcal{L}}$ where $\mathcal{K}=(K,
K^\circ)$ and $L_j = \{\emptyset, \{1\}, \{2\}\}$ for each $j\in
[m]$. }
\end{exam}

Definition~\ref{def:general-chessboard} can be specialized in many
ways. Again, we focus on the special cases motivated by intended
applications to the generalized Tverberg problem.

\begin{defin}\label{def:special-1}
Suppose that $\mathbf{k}=(k_i)_{i=1}^n$ and
$\mathbf{l}=(l_j)_{j=1}^m$ are two sequences of non-negative
integers. Then the complex,
\begin{equation}\label{eqn:defin-hv}
\Delta_{m,n}^{\mathbf{k}, \mathbf{l}} =
\Delta_{m,n}^{k_1,...,k_n;l_1,...,l_m}
\end{equation}
arises as the complex of all rook-placements $A\subset [m]\times
[n]$ such that at most $k_i$ rooks are allowed to be in the $i$-th
row (for $i=1,...,n$), and at most $l_j$ rooks are allowed to be
in the $j$-th column (for $j=1,...,m$).
\end{defin}

When $k_1=\cdots =k_n=p$ and $l_1=\cdots =l_m=q$, we obtain the
complex $\Delta_{m,n}^{p,q}$. For the reasons which will become
clear in the final section of the paper, we will be especially
interested in the case $l_1=\cdots =l_m=1$, i.e. in the complexes,
\begin{equation}\label{eqn:oznaka}
\Delta_{m,n}^{k_1,...,k_n;\mathbf{1}} :=
\Delta_{m,n}^{k_1,...,k_n;1,...,1}
\end{equation}

The inductive argument used in the proof of the main theorem
(Theorem~\ref{thm:main}) requires the analysis
(Proposition~\ref{prop:technical}) of complexes
$\Delta_{m,n}^{2,1(j)}$ which arise as follows. We assume that $R
\subset [n]$ is a $j$-element subset of $[n]$, prescribed in
advance, labelling selected $j$ rows in the $(m\times
n)$-chessboard.

\begin{defin}\label{def:special-2}
A rook-placement $A\subset [m]\times [n]$ is a simplex in
$\Delta_{m,n}^{2,1(j)}$ if and only if at most $2$ rooks are
allowed in rows indexed by $R$ and at most one in all other rows
and columns. Obviously for $j=0$ we obtain the usual chessboard
complex $\Delta_{m,n}^{2,1(0)}=\Delta_{m,n}$, and for $j=n$ the
generalized chessboard complex
$\Delta_{m,n}^{2,1(n)}=\Delta_{m,n}^{2,1}$. When $n=2$, aside from
these two possibilities, there is only one case remaining, the
complex $\Delta_{m,2}^{2,1(1)}$.
\end{defin}

\begin{rem}{\rm
The problem of determining the connectivity of generalized
chessboard complexes was considered in \cite{krw}, where they
proved a result implying that the homology with rational
coefficients $H_\nu(\Delta_{m,n}^{k_1,...,k_n;\mathbf{1}};
\mathbb{Q})$ is trivial if $\nu\leq (\mu -2)$ where
$$\mu =\min \{m,\left[\frac{m+k_1+\cdots
+k_n+1}3\right],k_1+\cdots +k_n\}.$$

If we are interested in the (homotopic) connectivity of
$\Delta_{m,n}^{k_1,...,k_n;\mathbf{1}}$ one can use the inductive
argument based on the application of the nerve lemma, used in
\cite{blvz}. By refining this argument we obtain here
(Theorem~\ref{thm:main}) a substantially better estimate in the
case $l_1=\cdots =l_m=1$. This is exactly the result needed here
for a proof of a generalized colored Tverberg theorem
(Theorem~\ref{thm:application}) for which the original estimate
from \cite{krw} was not sufficient. We believe and conjecture that
the same argument could be used to prove a better estimate in the
general (multidimensional) case. }
\end{rem}

For completeness and the reader's convenience here we state,
following \cite{bj}, a version of the {\em Nerve Lemma} needed in
the proof of the main theorem and other propositions.

\begin{lema}\label{lema:Nerve} {\bf (Nerve Lemma)}
\label{nerve1} Let $\Delta$ be a simplicial complex and $\{
L_i\}_{i=1}^k$ a family of subcomplexes such that $\Delta =
\cup_{i=1}^k~L_i.$ Suppose that every intersection $L_{i_1} \cap
L_{i_2} \cap \ldots \cap L_{i_t}$  is $(\mu-t+1)$-connected for $t
\geq 1.$ Then $\Delta$ is $\mu$-connected.
\end{lema}

\subsection{Selected examples of complexes $\Delta_{m,n}^{\mathcal{K}, \mathcal{L}}$}
\label{sec:examples}

As a preparation for the proof of Theorem~\ref{sec:main}, and as
an illustration of the use and versatility of the Nerve Lemma,
here we analyze in some detail the connectivity properties of
multiple chessboard complexes for some small values of $m$ and
$n$.

\begin{exam}
$\Delta_{m,1}^{2,1}=(\sigma^{m-1})^{(1)}$ is the $1$-skeleton of a
$(m-1)$-dimensional simplex, in particular it is connected.
\end{exam}

\begin{exam}
$\Delta_{3,2}^{2,1}\approx S^1\times I$, $\Delta_{4,2}^{2,1}$ has
the homology of $S^2$, $\Delta_{5,2}^{2,1}\approx S^3$, and for
$m\geq 5$ the complex $\Delta_{m,2}^{2,1}$ is $2$-connected.
\end{exam}

\medskip\noindent
{\bf Proof:} The complex $\Delta_{3,2}^{2,1}$ is a triangulation
of the surface of a cylinder into $6$ triangles.
$\Delta_{4,2}^{2,1}$ is a simplicial complex whose simplices are
subsets of the $4\times 2$ chessboard $\{(i,j) \mid 1\leq i\leq 4,
1\leq j\leq 2\}$ with at most two vertices in the same row and at
most one vertex in each column. This complex is covered by $4$
subcomplexes $L_1, L_2, L_3, L_4$ where $L_i$ is the collection of
simplices which contain $(i,1)$ as a vertex, together with their
faces. Each $L_i$ is contractible. For $i\neq j$, $L_i\cap L_j$ is
a union of a tetrahedron with two intervals joining the vertices
of the tetrahedron with two new vertices, hence it is also
contractible. The intersection of each three of these subcomplexes
is nonempty (the union of three intervals with a common vertex and
one additional vertex). Also, the intersection of all $4$
subcomplexes is a set of $4$ different points. Hence, by the Nerve
Lemma, the complex $\Delta_{4,2}^{2,1}$ is $1$-connected. Using
the Euler-Poincar\'{e} formula it is easy to see that this complex
has the homology of $S^2$.

We already know (Proposition~\ref{prop:Bier-1}) that
$\Delta_{5,2}^{2,1}$ is a $3$-sphere. However, as in the previous
example, this complex can be covered by $5$ contractible
subcomplexes $L_1,...,L_5$. The intersection of any two of these
subcomplexes is contractible. The intersection of any three of
them is also contractible (the union of three triangles with a
common edge and two additional edges connecting the vertices of
this edge with two additional points). The intersection of any
four of these subcomplexes as well as the intersection of all of
them is non-empty. Therefore, by Nerve lemma, the complex
$\Delta_{5,2}^{2,1}$ is $2$-connected.

It is easy to verify that this complex is a simplicial
$3$-manifold; the link of each vertex is a $2$-dimensional sphere,
the link of each edge is a circle, and the link of each
$2$-dimensional simplex is a $0$-dimensional sphere. Hence,
$\Delta_{5,2}^{2,1}\approx S^3$. \hfill $\square$

\begin{exam}
$\Delta_{3,3}^{2,1}$ and $\Delta_{4,3}^{2,1}$ are $1$-connected,
and $\Delta_{5,3}^{2,1}$ and $\Delta_{6,3}^{2,1}$ are
$2$-connected.
\end{exam}

\medskip\noindent
{\bf Proof:} It is easy to see that each simplex in these
complexes is a face of a simplex having a vertex in the first
column. As before we apply the Nerve Lemma to the covering by
subcomplexes $M_i$, $i\in \{1,2,3\}$, where $M_i$ is the
collection of simplices having vertex at $(1,i)$, together with
their faces. \hfill $\square$
\bigskip

\begin{exam}
The complex $\Delta_{m,2}^{2,1(1)}$
(Definition~\ref{def:special-2}) is connected for $m\geq 2$, and
$1$-connected for $m\geq 4$.
\end{exam}

\medskip\noindent
{\bf Proof:} The proof goes along the same lines as before and
uses the Nerve lemma. As we already know
(Proposition~\ref{prop:Bier-1}) $\Delta_{4,2}^{2,1(1)}$ is a
simplicial surface homeomorphic to $S^2$. This can be proved
directly as follows. The link of each vertex is a circle (a
triangle or a hexagon), and the link of each edge consists of two
points, i.e. it is $S^0$. The Euler characteristic of this complex
equals $2$, so it must be $S^2$. \hfill $\square$

\bigskip
Some of the examples of generalized chessboard complexes are
spheres (Proposition~\ref{prop:Bier-1} and
Example~\ref{exam:Bier}). Here we meet some quasi-manifolds.

\begin{exam}
$\Delta_{7,3}^{2,1}$ is a $5$-dimensional quasi-manifold. More
generally $\Delta_{pn+1,n}^{p,1}$ is a $(pn-1)$-dimensional
quasi-manifold for each $p\geq 1$.
\end{exam}

\medskip\noindent
{\bf Proof:} It is easy to verify that the link of each
$4$-dimensional simplex is $S^0$, the link of each $3$-dimensional
simplex is a combinatorial circle (consisting of either $3$ or $6$
edges), the link of some $2$-dimensional simplices is $S^2$
(triangles with two vertices in the same row), but some other
$2$-dimensional simplices (whose vertices are in three different
rows) have the link homeomorphic to the torus rather than to
$S^2$. A similar proof applies in the general case.
\hfill$\square$
\bigskip

Concerning the generalized chessboard complexes
$\Delta_{m,n}^{p,q}$ with higher values of $q$, we mention one
additional simple example.
\bigskip

\begin{exam}
$\Delta_{m,n}^{m-1,n}\approx S^{(m-1)n-1}$, and
$\Delta_{m,n}^{p,n}$ is $(np-2)$-connected.
\end{exam}

\medskip\noindent
{\bf Proof:} $\Delta_{m,1}^{m-1,1}$ is the boundary of
$(m-1)$-dimensional simplex, and so homeomorphic to $S^{m-2}$.
Since $\Delta_{m,n}^{m-1,n}$ is a join of $n$ copies of
$\Delta_{m,1}^{m-1,1}$, it is homeomorphic to $S^{(m-1)n-1}$.

Similarly we see that the complex $\Delta_{m,n}^{p,n}$ is a join
of $n$ complexes of the type $\Delta_{m,1}^{p,1}$ which is
identified as the $(p-1)$-skeleton of an $(m-1)$-dimensional
simplex. We conclude that this complex is a wedge of
$(np-1)$-dimensional spheres, so must be $(np-2)$-connected. In
particular for $p=1$ this reduces to the fact that
$\Delta_{m,n}^{1,n}$ is a join of $n$ copies of the finite set of
$m$ points, and so $(n-2)$-connected. \hfill $\square$
\bigskip

\begin{cor} The complex $\Delta_{3,2}^{2,2}$ is a $3$-sphere,
 $\Delta_{3,2}^{2,2}\cong S^3$.
\end{cor}
\section{Connectivity of multiple chessboard complexes}
\label{sec:main}

Theorem~\ref{thm:main}, as the first of the two main result of our
paper, provides an estimate for the connectivity of the
generalized chessboard complex
$\Delta_{m,n}^{k_1,...,k_n;\mathbf{1}}$
(Definition~\ref{def:special-1} and equation (\ref{eqn:oznaka})).

\medskip
By the Hurewicz theorem in order to show that a connected complex
$K$ is $k$-connected ($k\geq 1$) it is sufficient to show that
$\pi_1(K)=0$ and that $\widetilde{H}_j(K; \mathbb{Z}) = 0$ for
each $j=1,\ldots, k$.

\begin{prop}\label{prop:main-intro}
If both $m\geq n+2$ and $k_1+\ldots + k_n \geq 3$ then,
\begin{equation}\label{eqn:fundam}
\pi_1(\Delta_{m,n}^{k_1,...,k_n;\mathbf{1}}) = 0.
\end{equation}
\end{prop}

\medskip\noindent
{\bf Proof:} If $k_1=\ldots = k_n = 1$ (the case of the standard
chessboard complex) the condition (\ref{eqn:fundam}) reduces to
$m-2 \geq n \geq 3$ which is (following \cite{blvz}) sufficient
for the $1$-connectivity of $\Delta_{m, n}^{1,...,1;
\mathbf{1}}\cong \Delta_{m, n}$. Small examples of generalized
chessboard complexes, as exhibited in Section~\ref{sec:examples},
also support the claim in the case $n=2$.

\medskip The general case is established by the {\em Gluing Lemma}
\cite{bj} (or Seifert-van Kampen theorem) following the ideas of
similar proofs \cite[Theorem~1.1]{blvz}  and \cite[Theorem
3]{zv92}. Recall that the Gluing Lemma is essentially the $k=2$
case of the Nerve Lemma (Lemma~\ref{lema:Nerve}). \hfill $\square$

\begin{theo}\label{thm:main}
\label{prva} The generalized chessboard complex
$\Delta_{m,n}^{k_1,...,k_n;\mathbf{1}}$ is $(\mu-2)$-connected
where,
$$ \mu = \min \{m-n+1,k_1+\cdots +k_n\}.$$ In particular if $m\geq
k_1+\cdots +k_n+n-1$ then $\Delta_{m,n}^{k_1,...,k_n;\mathbf{1}}$
is $(k_1+\cdots +k_n-2)$-connected.
\end{theo}

\medskip\noindent
{\bf Proof:} By the Hurewicz theorem and
Proposition~\ref{prop:main-intro} it is sufficient to show  the
complex $\Delta_{m,n}^{k_1,...,k_n;\mathbf{1}}$ is homologically
$(\mu-2)$-connected.

We carry on the proof by showing that the complex
$\Delta_{k_1+\cdots +k_n+n-1,n}^{k_1,...,k_n;\mathbf{1}}$ is
$(k_1+\cdots +k_n-2)$-connected, and that by reducing the number
of columns by $1$ the connectivity degree of the complex either
reduces by $1$ or remains the same.

We proceed by induction. It is easy to check that the statement of
the theorem is true for $n=1$ and every $m$ and $k_1$, and that
our estimate is true if $m\leq n$. It also follows directly from
the known result from \cite{zv92} when $k_1=\cdots =k_n=1$. Let us
suppose that the statement is true for the complex
$\Delta_{r,s}^{k_1,...,k_s;1}$ when $s<n$ (for every $r$ and
$k_1,...,k_s$), and also for $s=n$ if $r<m$.

\medskip
We now focus attention to the complex $\Delta =
\Delta_{m,n}^{k_1,...,k_n;\mathbf{1}}$, taking into the account
that the case when all numbers $k_1,...,k_n$ are equal to $1$ is
already covered.

\medskip

{\bf (i)} Let us start with the case $m\geq k_1+\cdots +k_n+n-1$.
Note that in this case $\mu = k_1+\cdots +k_n$. Without loss of
generality (by permuting the rows if necessary) we can assume that
$k_1\geq \ldots\geq k_n$, in particular we can assume that
$k_1\geq 2$.

\medskip
The complex $\Delta = \Delta_{m,n}^{k_1,...,k_n;\mathbf{1}}$ is
covered by the contractible subcomplexes (cones) $L_i = {\rm
Star}_{\Delta}(v_i)$ having the apices at the points $v_i=(i,1)$,
$i=1,...,m$. Let us analyze the connectivity degree of
intersections $L_i\cap L_j$ (where without loss of generality we
assume that $i = m-1$ and $j=m$).

\medskip  The intersection $L_{m-1}\cap L_m = K \cup K'$ is
the union of two sub-complexes where,
\[
K = \{A\in \Delta \mid A\cup \{(m-1,1), (m,1)\}\in \Delta\}
\]
and $B\in K'$ if and only if for some $A\supseteq B$,
\begin{equation}\label{eqn:K-prim}
A\cup \{(m-1,1)\}\in \Delta,\quad  A\cup \{(m,1)\}\in \Delta,\quad
A\cup \{(m-1,1), (m,1)\}\notin \Delta.
\end{equation}
The last condition in (\ref{eqn:K-prim}) can be replaced by the
condition that $\vert A\cap ([m-2]\times\{1\}) \vert = k_1-1$
i.e.\ that $A\subset [m-2]\times [n]$ has $k_1-1$ elements in the
first row.

\medskip
By construction $K \cong \Delta_{m-2,n}^{k_1-2,...,k_n;\mathbf{1}}
\ast I^1$. The complex $K'$ is a subcomplex of $\Delta' =
\Delta_{m-2,n}^{k_1-2,...,k_n;\mathbf{1}}$ (based on the
chessboard $[m-2]\times [n]$) and its structure is described by
the following lemma.

\medskip\noindent
{\bf Lemma:} Let $\mathcal{T}$ be the collection of all
$(k_1-1)$-element subsets of $[m-2]\times\{1\}$ and for $T\in
\mathcal{T}$ let $K_T = {\rm Star_{\Delta'}(T)}$ where $\Delta' =
\Delta_{m-2,n}^{k_1-2,...,k_n;\mathbf{1}}$. Then $K' = \cup_{T\in
\mathcal{T}}~K_T$. Moreover, for any proper subset
$\mathcal{T}'\subset \mathcal{T}$ and $S\in \mathcal{T}\setminus
\mathcal{T}'$
\begin{equation}\label{eqn:proper}
(K \cup \bigcup_{T\in \mathcal{T}'}~K_T)\cap K_S = K\cap K_S.
\end{equation}

\medskip
We continue the analysis of the complex $K'$ by observing that the
complex $K\cap K_S$ is the join of the boundary of the simplex $S$
(homeomorphic to the sphere $S^{k_1-3}$) and a complex isomorphic
to $\Delta_{m-k_1-1,n-1}^{k_2,...,k_n;\mathbf{1}}$ (which is
$(k_2+\cdots +k_n-2)$-connected by the induction hypothesis).
Hence, this complex is $(\mu -4)$-connected.

\medskip
By Lemma the complex $K'$ can be built from $K$ by adding
complexes $K_S$, one at a time, where $S\in \mathcal{T}$. By using
repeatedly the Mayer-Vietoris long exact sequence (or
alternatively the Gluing Lemma), we see that the complex
$L_{m-1}\cap L_m$ is $(\mu -3)$-connected.

\medskip
For the reader's convenience, before we proceed to the general
case, we prove that the intersection of any three of the
subcomplexes $L_i$, $i\in \{1,...,m\}$ is $(\mu -4)$-connected.

\medskip
We begin by analyzing how the maximal simplices $A$ in
$L_{m-2}\cap L_{m-1}\cap L_m$ are allowed to look like. More
precisely we look at the intersection $A\cap ([m]\times \{1\})$ of
$A$ with the first row of the chessboard and classify $A$
depending on the size of the sets $A\cap ([m-3]\times \{1\})$ and
$A\cap (\{m-2, m-1, m\}\times \{1\})$ (lets denote these
cardinalities by $\alpha_A$ and $\beta_A$ respectively). We
observe that the case $\alpha_A = k_1-3$ and $\beta_A = 3$ is the
first (and only) case when $\alpha_A+\beta_A =k_1$. Indeed, if
$\beta_A< 3$ (say for example $A\cap (\{m-2, m-1, m\}\times \{1\})
= (\{m-2, m-1\}\times \{1\})$) then $\alpha_A + \beta_A<k_1$,
since otherwise $A\notin L_m$. From here it immediately follows
that the case $\beta_A=2$ can be excluded being subsumed by the
case $\beta_A = 3$.

\bigskip
Summarizing we observe that  $L_{m-2}\cap L_{m-1}\cap L_m$ can be
built from the complexes $K$, $K'$ and $K''$ (if $k_1\geq 3$)
generated (respectively) by simplices $A$ with $(\alpha_A,
\beta_A) = (k_1-3,3), (\alpha_A, \beta_A) = (k_1-2,1)$ or
$(\alpha_A, \beta_A) =(k_1-1,0)$.

\medskip
The complex $K$ is a join of the triangle (with vertices
$(m-2,1),(m-1,1),(m,1)$) and the complex of the type
$\Delta_{m-3,n}^{k_1-3,...,k_n;\mathbf{1}}$, hence it is
contractible. The complex $K'$ is a join of a finite set of $3$
points ($(m-2,1),(m-1,1),(m,1)$) and the complex of the type
$\Delta_{m-3,n}^{k_1-2,...,k_n;\mathbf{1}}$. Let $\mathcal{T}$ be
the collection of all $(k_1-2)$-element subsets of $[m-3]$. Then
$K'$ can be represented as the union $K' = \cup_{T\in
\mathcal{T}}~K'_T$ where,
\begin{equation}\label{eqn:3-inter}
K'_T \cong \Delta_{m-k_1-1,n-1}^{k_2,...,k_n;\mathbf{1}} \ast
\Delta[T]\ast [3]
\end{equation}
($\Delta[S]$ is the simplex spanned by vertices in $S$). In words,
the complex $K'$ can be represented as the union of the complexes
which are joins of the three point set with the simplex
$\Delta^{k_1-3}$ of dimension $k_1-3$ and with the complex of the
type $\Delta_{m-k_1-1,n-1}^{k_2,...,k_n;\mathbf{1}}$.

\medskip
The complexes $K'_T$ are contractible and we observe that by
adding them to $K$, one by one, the intersection of each of them
with the previously built complex is of the type $[3]\ast
S^{k_1-4}\ast \Delta_{m-k_1-1,n-1}^{k_2,...,k_n;\mathbf{1}}$. This
complex is $(\mu -4)$-connected, consequently $K\cup K'$ is
$(\mu-4)$-connected as well.

\medskip
The complex $K''$ is the complex of the type
$\Delta_{m-3,n}^{k_1-1,...,k_n;\mathbf{1}}$ and can be represented
as the union of the complexes $K''_T$ which are joins of the
simplex $\Delta[T]\cong \Delta^{k_1-2}$ with the complex of the
type $\Delta_{m-k_1-2,n-1}^{k_2,...,k_n;\mathbf{1}}$. These
complexes are contractible and when adding them to  $K\cup K'$,
one by one, the intersection is of the type $S^{k_1-3}\ast
\Delta_{m-k_1-2,n-1}^{k_2,...,k_n;\mathbf{1}}$. This complex is
$(\mu -5)$-connected.

By using repeatedly the Mayer-Vietoris long exact sequence we
finally observe that the complex $L_{m-2}\cap L_{m-1}\cap L_m =
K\cup K'\cup K''$ is $(\mu -4)$-connected.

\medskip
Now we turn to the general case. We want to prove that the
intersection of any collection of $q$ cones $L_i$, let us say
$L_{m-q+1}\cap L_{m-q+2}\cap \cdots \cap L_m$, is $(\mu
-q-1)$-connected. We treat separately the cases $q\leq k_1$, and
$q>k_1$. (As before we are allowed to assume that $k_1\geq 2$.)

{\it (a)}\quad $q\leq k_1$: This case is treated similarly as the
special cases $q=2$ and $q=3$. Our objective is to prove that the
intersection $L^{(q)} = L_{m-q+1}\cap L_{m-q+2}\cap \cdots \cap
L_m$ is $(\mu -q-1)$-connected.

\medskip\noindent
We begin again by analyzing how the maximal simplices $A$ in
$L^{(q)}$ are allowed to look like by looking at the pairs of
integers $(\alpha_A, \beta_A)$ where,
 $$
 \alpha_A =\vert A\cap ([m-q]\times \{1\})\vert \quad \mbox{\rm
 and} \quad \beta_A = \vert A\cap (\{m-q+1, \ldots m\}\times
\{1\})\vert.
 $$
As before we observe that (for maximal $A$) the case $\alpha_A +
\beta_A = k_1$ is possible only if  $\beta_A = q$ and $\alpha_A =
k_1-q$. Moreover we observe that the intersection $L^{(q)}$ can be
expressed as the union of the complexes $K_1,...,K_q$ where $K_1$
is generated by simplices $A$ of he type $(\alpha_A, \beta_A) =
(k_1-q, q)$ while for $j>1$ the complex $K_j$ is generated by
simplices of the type $(\alpha_A, \beta_A) = (k_1-q+j-1, q-j)$.

\medskip\noindent
The complex $K_1$ is the join of the $(q-1)$-dimensional simplex
and the complex of the type
$\Delta_{m-q,n}^{k_1-q,k_2,...,k_n;\mathbf{1}}$ (or of the type
$\Delta_{m-q,n-1}^{k_2,...,k_n;\mathbf{1}}$ if $q=k_1$). So, $K_1$
is contractible.

\medskip\noindent
The complex $K_2$ can be presented as the union of complexes which
are joins of the $(q-3)$-dimensional skeleton of the
$(q-1)$-dimensional simplex with the $(k_1-q)$-dimensional
simplex, and with the complex of the type
$\Delta_{m-k_1-1,n-1}^{k_2,...,k_n;\mathbf{1}}$. These complexes
are contractible, and when adding one by one to the complex $K_1$
we notice that the intersection of each of them with the
previously built complex is the complex of the type of the join of
the $(q-3)$-dimensional skeleton of the $(q-1)$-dimensional
simplex with the $(k_1-q-1)$-dimensional skeleton of the
$(k_1-q)$-dimensional simplex, and with the complex of the type
$\Delta_{m-k_1-1,n-1}^{k_2,...,k_n;\mathbf{1}}$. This intersection
is, by the induction hypothesis, $(\mu -4)$-connected. So, the
union $K_1\cup K_2$ is $(\mu -3)$-connected.

\medskip\noindent
We proceed in the same way by adding complexes $K_j$, one at the
time. Finally, $K_q$ is the subcomplex consisting of simplices
having no vertices in the set $\{(m-q+1,1),...,(m,1)\}$, and so it
is of the type $\Delta_{m-q,n}^{k_1-1,...,k_n;\mathbf{1}}$. The
complex $K_q$ could be presented as the union of complexes which
are joins of the $(k_1-2)$-dimensional simplex with the complex of
the type $\Delta_{m-q-k_1+1,n-1}^{k_2,...,k_n;\mathbf{1}}$. These
complexes are contractible, and when adding one by one to the
complex $K_1\cup ...\cup K_{q-1}$ we notice that the intersection
of each of them with the previously built complex is the complex
of the type of the join of the $(k_1-3)$-dimensional skeleton of
the $(k_1-2)$-dimensional simplex with the complex of the type
$\Delta_{m-q-k_1+1,n-1}^{k_2,...,k_n;\mathbf{1}}$. This
intersection is, by the induction hypothesis,
$(\mu-q-2)$-connected. So, the union $K_1\cup ...\cup K_{k_1}$ is
$(\mu-q-1)$-connected.

\medskip
{\it (b)}\quad $q>k_1$: Let us prove that the intersection of any
$q$ cones, for example the intersection $L_{m-q+1}\cap
L_{m-q+2}\cap \cdots \cap L_m$, is $(\mu -q-1)$-connected in this
case as well. As before we express this intersection as the union
of  complexes $K_1,...,K_{k_1}$. Here, $K_1$ is the subcomplex
consisting of simplices having $k_1-1$ vertices in the set
$\{(m-q+1,1),...,(m,1)\}$ and it is the complex of the type of
join of the $(k_1-2)$-dimensional skeleton of the
$(q-1)$-dimensional simplex and the complex of the type
$\Delta_{m-q,n-1}^{k_2,...,k_n;\mathbf{1}}$. So, it is $(\mu
-q+k_1-2)$-connected by the induction hypothesis. The complex
$K_2$ consists of simplices having $k_1-2$ vertices in the set
$\{(m-q+1,1),...,(m,1)\}$ and one vertex in the remaining vertices
of the first row. The type of this complex is the join of the
$(k_1-3)$-dimensional skeleton of the $(q-1)$-dimensional simplex,
and the complex of the type
$\Delta_{m-q,n}^{1,k_2,...,k_n;\mathbf{1}}$. The complex $K_2$ can
be presented as the union of complexes which are joins of the
$(k_1-3)$-dimensional skeleton of the $(q-1)$-dimensional simplex
with a point, and with the complex of the type
$\Delta_{m-q-1,n-1}^{k_2,...,k_n;\mathbf{1}}$. These complexes are
contractible, and when adding one by one to the complex $K_1$ we
notice that the intersection of each of them with the previously
built complex is the complex of the type of the join of the
$(k_1-3)$-dimensional skeleton of the $(q-1)$-dimensional simplex,
and the complex of the type
$\Delta_{m-q-1,n-1}^{k_2,...,k_n;\mathbf{1}}$. This intersection
is, by the induction hypothesis, $(\mu-q+k_1-4)$-connected. So,
the union $K_1\cup K_2$ is $(\mu-q+k_1-3)$-connected.

We proceed in the same way. Finally, $K_{k_1}$ is the subcomplex
consisting of simplices having no vertices in the set
$\{(m-q+1,1),...,(m,1)\}$, and so it is of the type
$\Delta_{m-q,n}^{k_1-1,...,k_n;\mathbf{1}}$. The complex $K_{k_1}$
can be presented as the union of complexes which are joins of the
$(k_1-2)$-dimensional simplex with the complex of the type
$\Delta_{m-q-k_1+1,n-1}^{k_2,...,k_n;\mathbf{1}}$. These complexes
are contractible, and when adding one by one to the complex
$K_1\cup ...\cup K_{k_1-1}$ we notice that the intersection of
each of them with the previously built complex is the complex of
the type of the join of the $(k_1-3)$-dimensional skeleton of the
$(k_1-2)$-dimensional simplex with the complex of the type
$\Delta_{m-q-k_1+1,n-1}^{k_2,...,k_n;\mathbf{1}}$. This
intersection is, by the induction hypothesis,
$(\mu-q-2)$-connected. So, the union $K_1\cup ...\cup K_{k_1}$ is
$(\mu-q-1)$-connected.
\smallskip

{\bf (ii)} Consider now the case $k_1+\cdots +k_n\leq m < 
k_1+\cdots +k_n+n-1$ and let us prove that the complex is
$(m-n-1)$-connected. We again cover the complex by the cones
$L_i$, $i=1,...,m$, and note that the intersection of any two of
them (let us say $L_{m-1}$ and $L_m$) could be built by adding
contractible subcomplexes so that the intersection of any of them
with previously built complex is of the type $S^{k_1-3} \ast
\Delta_{m-k_1-1,n-1}^{k_2,...,k_n;\mathbf{1}}$. This complex is
$(m-n-3)$-connected by the induction hypothesis, and so the
intersection $L_{m-1}\cap L_m$ is $(m-n-2)$-connected by the
Mayer-Vietoris theorem. In the same way we prove that the
intersection of any three cones is $(m-n-3)$-connected etc.
\smallskip

{\bf (iii)} In the case $m\leq k_1+\cdots +k_n-1$ we cover the
complex $\Delta_{m,n}^{k_1,...,k_n;\mathbf{1}}$ by the cones $M_i$
with the vertices $(1,i)$. The intersection of any two of them
(let us say $M_{n-1}$ and $M_n$) is the complex of the type
$\Delta_{m-1,n}^{k_1,...,k_{n-1}-1,k_n-1;1}$ (if $k_{n-1},k_n\geq
2$), or $\Delta_{m-1,n-1}^{k_1,...,k_{n-1}-1;1}$ (if exactly one
of them, let us say $k_n$, equals $1$), or
$\Delta_{m-1,n-2}^{k_1,...,k_{n-2};1}$ (if $k_{n-1}=k_n=1$). In
any case, this intersection is ar least $(m-n-2)$-connected by the
induction hypothesis. In the same way we prove that the
intersection of any three cones is $(m-n-3)$-connected etc. \hfill
$\square$

\begin{cor}
\label{druga'} $\Delta_{m,n}^{2,1}$ is $(2n -2)$-connected for
$m\geq 3n-1$, $\Delta_{m,n}^{3,1}$ is $(3n -2)$-connected for
$m\geq 4n-1$, and generally $\Delta_{m,n}^{p,1}$ is $(pn
-2)$-connected for $m\geq (p+1)n-1$.
\end{cor}

Notice that the general estimate obtained in \cite{krw} implies
that the complex $\Delta_{m,n}^{p,1}$ is $(pn-2)$-connected for
$m\geq 2pn-1$, which is (compared to $m\geq (p+1)n-1$) a weaker
estimate (roughly by a factor of $2$).

\begin{cor}
$\Delta_{7,3}^{2,1}$ is $3$-connected, but not $4$-connected.
\end{cor}

\medskip\noindent
{\bf Proof:} For the proof of the last claim it suffices to
compute the Euler characteristic of this complex $\chi
(\Delta_{7,3}^{2,1})=147$. Since $\Delta_{7,3}^{2,1}$ is
$3$-connected $5$-dimensional quasi-manifold, we have $\beta_5=1$
and so $\beta_4=147$. \hfill $\square$

\begin{rem}{\rm
The estimate $m-n-1$ for small values of $m$ in the statement of
Theorem \ref{prva} can be significantly improved. For example, the
following result gives the estimate on the connectivity of the
generalized chessboard complex
$\Delta_{m,n}^{k_1,...,k_n;\mathbf{1}}$ in the case $k_1=\cdots
=k_j=2$ and $k_{j+1}=\cdots =k_n=1$, i.e. when $j$ of $n$ numbers
$k_1,...,k_n$ are equal to $2$ and the remaining $n-j$ are equal
to $1$. We believe that this estimate is close to the best
possible. Recall (Definition~\ref{def:special-2}) that this
complex is already introduced as the complex
$\Delta_{m,n}^{2,1(j)}$.}
\end{rem}

\begin{prop}\label{prop:technical}
The complex $\Delta_{m,n}^{2,1(j)}$
(Definition~\ref{def:special-2}) is $(\mu -2)$-connected where

$$\mu=\left\{\begin{array}{rl}
m, & m < \frac{n+j}2\\
\left[ \frac{m+n+j+1}3\right], & \frac{n+j}2\leq m< n+\frac j2-1\\
\left[ \frac{5m+n+2j+5}9\right], & n+\frac j2 -1\leq m < n+2j\\
\left[ \frac{m+n+2j+1}3\right], & n+2j\leq m < 2n+j-1\\
n+j, & m\geq 2n+j-1
\end{array}
\right. $$
\end{prop}

The proof uses exactly the same ideas, so we omit the details.

As a final comment we repeat that, motivated by possible
applications to theorems of Tverberg type, we are interested
mostly in the values of $m$ for which the complex
$\Delta_{m,n}^{k_1,...,k_n;\mathbf{1}}$ is $(k_1+\cdots
+k_n-2)$-connected. We believe that the assumption $m\geq
k_1+\cdots +k_n+n-1$ is optimal in that respect.

\section{Shellability of multiple chessboard complexes}
\label{sec:shelling-short}

For the definition and basic facts about shellable complexes the
reader is referred to \cite{Jo-book} and \cite{Koz}. One of the
central topological properties of these complexes is the following
well known lemma.

\begin{lema}\label{lema:shell}
A shellable, $d$-dimensional simplicial complex is either
contractible or homotopy equivalent to a wedge of $d$-dimensional
spheres.
\end{lema}
An immediate consequence of Lemma~\ref{lema:shell} is that a
$d$-dimensional, shellable complex is always $(d-1)$-connected.
This observation opens a way of proving Theorem~\ref{thm:main} by
showing that the associated multiple chessboard complex is
shellable.

\medskip
Shellability of standard chessboard complexes $\Delta_{m,n}$ for
$m\geqslant 2n-1$ is established by G. Ziegler in \cite{zie}. He
established {\em vertex decomposability} of these and related
complexes, emphasizing that the natural lexicographic order of
facets of $\Delta_{m,n}$ is NOT a shelling order. We demonstrate
that a version of {\em `cyclic reversed lexicographical order'} is
a shelling order both for standard and for generalized chessboard
complexes. Before we prove the general case
(Theorem~\ref{thm:shelling-main}), we outline the main idea by
describing a shelling order for the standard chessboard complexes
$\Delta_{m,n}$.

\bigskip\noindent {\bf Shelling order for $\Delta_{m,n}$:} Let $\Delta_{m,n}$,
be a chessboard complex which satisfies the condition $m\geq
2n-1$. If $A = (a_1, a_2,\ldots, a_n )$ is a sequence of distinct
elements of $[m]$ the associated simplex $\{(a_i, i)\}_{i=1}^n$ in
$\Delta_{m,n}$ is denoted by $\hat{A}$. Both $A$ and $\hat{A}$ are
interchangeably referred to as facets of $\Delta_{m,n}$.

\medskip
The shelling order $\ll$ on $\Delta_{m,n}$ is introduced by
describing a rule (algorithm) which decides for each two distinct
facets $A$ and $B$ of $\Delta_{m,n}$ whether $A\ll B$ or $B\ll A$.
We adopt a basic cyclic order $\prec$ on $[m]$,
\[
1 \prec 2 \prec \ldots \prec m \prec 1
\]
which for each $a\in [m]$ reduces to a genuine linear order
$\prec_a$ on [m],
\begin{equation}\label{eqn:stand-linar}
a+1 \prec_a a+2 \prec_a \ldots \prec_a m \prec_a 1 \prec_a \ldots
\prec_a a,
\end{equation}
and in particular $\prec_m$ is the standard linear order $<$ on
$[m]$.

\medskip
Suppose that $A = (a_1,\ldots, a_n)$ and $B = (b_1,\ldots, b_n)$
are two distinct facets of $\Delta_{m,n}$. The procedure of
comparing $A$ and $B$ begins by comparing $a_1$ and $b_1$. By
definition the relation $A\ll B$ is automatically satisfied if
$a_1<b_1$. If $a_1=b_1 = a$ we use the order $\prec_{a}$ to
compare $\hat{A}$ and $\hat{B}$, first in the column $a-1$ then
(if necessary) in column $a-2$, then (if necessary) in column
$a-3$, etc. More precisely let us define the `comparison interval'
$[a-(p-1),a-1 ]$ by the requirement that,
\begin{enumerate}
 \item[(1)] for each $b\in [a-(p-1),a-1 ]$ both $\hat{A}$ and
$\hat{B}$ have a rook in the column $b$;
 \item[(2)] either $\hat{A}$ or $\hat{B}$ (or both) have no rooks
 in the column $a-p$.
\end{enumerate}

\medskip
We compare, moving from right to left (descending in the order
$\prec_{a}$), the positions of rooks of facets $\hat{A}$ and
$\hat{B}$ in the smaller chessboard $[a-(p-1), a-1]\times [n]$. If
$\{e\}\times [n]$ is the first column where they disagree, say
$\hat{A}\ni (e,i)\neq (e,j)\in \hat{B}$, then $A\ll B$ if $i>j$
(the rook corresponding to $(e,i)$ is above the rook associated to
$(e,j)$).

\medskip
Alternatively the facets $\hat{A}$ and $\hat{B}$ agree on the
whole of the smaller chessboard $[a-(p-1), a-1]\times [n]$. In
this case by definition $A\ll B$ if $\hat{B}$ does have a rook in
column $a-p$ and $\hat{A}$ does not.

\medskip
The final possibility is that the facets $\hat{A}$ and $\hat{B}$
agree on the comparison interval and neither $\hat{A}$ nor
$\hat{B}$ have a rook in the column $a-p$. If this is the case we
declare that the {\em first stage} of the comparison procedure is
over and pass to the second stage.

\medskip
The second stage of the comparison procedure begins by removing
(or simply ignoring) the first row of the chessboard $[m]\times
[n]$ and the column $\{a-p\}\times [n]$, together with the small
chessboard $[a-(p-1), a-1]\times [n]$ and all the rows of
$[m]\times [n]$ associated to the rooks in,
\[
([a-(p-1), a-1]\times [n])\cap \hat{A} = ([a-(p-1), a-1]\times
[n])\cap \hat{B}.
\]
This way we obtain a new chessboard $M\times N\subset [m]\times
[n]$ which inherits the (cyclic) `right to the left' order in each
of the rows so we can continue by applying the first stage of the
comparison procedure to the facets $\hat{A}' = \hat{A}\cap
(M\times N)$ and $\hat{B}' = \hat{B}\cap (M\times N)$.

\medskip
This process can be continued, stage after stage, and if $A\neq B$
sooner or later will lead to the decision whether $A\ll B$ or
$B\ll A$.

\begin{lema}\label{lema:linear}
The relation $\ll$ is a linear order on the set of facets of the
chessboard complex $\Delta_{m,n}$.
\end{lema}

\medskip\noindent
{\bf Proof of Lemma~\ref{lema:linear}:} By construction  $A\ll B$
and $B \ll A$ cannot hold simultaneously so it is sufficient to
show that the relation $\ll$ is transitive. Assume that $A\ll B$
and $B\ll C$. If both inequalities are decided at the same stage
then by inspection of the priorities one easily deduces that $A\ll
C$. Suppose that the inequalities are decided at different stages,
say $A\ll B$ is decided at the stage $i$ and $B\ll C$ at the stage
$j$ where (for example) $i<j$. Since $B$ and $C$ are not
discernible from each other, up to the stage $i$, we conclude that
the same argument used to decide that $A\ll B$ leads to the
inequality $A\ll C$.
 \hfill $\square$

\medskip
The details of the proof that the described linear order of facets
of $\Delta_{m,n}$ is indeed a shelling are omitted since they will
appear in greater generality in the proof of
Theorem~\ref{thm:shelling-main}. Nonetheless, in the following
very simple example we pinpoint the main difference between this
linear order and the lexicographic order of facets.

\begin{exam}\label{exam:sirina-2}
The lexicographic order of facets in the chessboard complex
$\Delta_{m,2}$ (for $m\geq 3$) is not a shelling. Indeed, if $A$
is a predecessor of the simplex $B = \{(2,1),(1,2)\}$ in the
lexicographic order then $A\cap B = \emptyset$. On the other hand
in the shelling order described above each of the simplices
$\{(2,1),(j,2)\}$ for $j\geq 3$ is a predecessor of $\sigma$.
\end{exam}

\medskip
\begin{theo}\label{thm:shelling-main}
For $m\geq k_1+k_2+\cdots+k_n+n-1$ the complex
$\Delta^{k_1,\ldots,k_n;\mathbf{1}}_{m,n}$ is shellable.
\end{theo}

\medskip\noindent
{\bf Proof:} Facets of $\Delta^{k_1,\ldots,k_n;\mathbf{1}}_{m,n}$
are encoded as $n$-tuples $A=(A_1,A_2,\ldots,A_n)$ of disjoint
subsets of $[m]$ where $|A_i|=k_i$. The elements of $A_i$
represent the positions of rooks in the $i$-th row, so strictly
speaking the simplex associated to $A$ is the set $\hat{A} = \{(a,
i)\in [m]\times [n] \mid a\in A_i\}$.  Both $A$ and $\hat{A}$
unambiguously refer to the same facet of $\Delta_{m,n}$.

\medskip
The linear order $\ll$ of facets of the multiple chessboard
complex $\Delta^{k_1,\ldots,k_n;\mathbf{1}}_{m,n}$ is defined by a
recursive procedure which generalizes the procedure already
described in the case of standard chessboard complex
$\Delta_{m,n}$.

\medskip
If $A_1$ is anti-lexicographically less than $B_1$ in the sense
that $max (A_1\triangle B_1)\in B_1$ we declare that $A\ll B$. If
$A_1=B_1=\{a_1,a_2,\ldots,a_{k_1}\}$ we consider $[m]\setminus
A_1=I_1\cup I_2 \cup\cdots\cup I_r$ where,
$$I_j=\{x_j,x_j+1,\cdots,x_j+s_j\}\subset [m]$$
are maximal sets (lacunas) of consecutive integers in
$[m]\setminus A_1$. We assume that $\textrm{max }I_j< \textrm{min
} I_{j+1}$ for all $j=1,2,\ldots,r-1$.

In full agreement with (\ref{eqn:stand-linar}) we introduce a
priority order $\prec$ of elements contained in the union of all
lacunas $I_j$ as follows. The priority order within the lacuna
$I_j$ is from `right to left',
\begin{equation}\label{E:lacuna}
x_j+s_j \succ x_j+s_j-1 \succ \ldots\succ  x_j+1 \succ x_j.
\end{equation}
The priority order of lacunas themselves is from `left to right',
so summarizing, the elements of $\cup_{j=1}^r~I_j$ listed in the
priority order $\succ$ from the biggest to the smallest are the
following,
\begin{equation}\label{E:order}
x_1+s_1,x_1+s_1-1,\ldots,x_1+1,x_1,x_2+s_2,\ldots,x_2,\ldots,x_r+s_r,\ldots,x_r.
\end{equation}
\\
If $A_1=B_1$ we define $A\ll B$ if either of the following
conditions is satisfied.
\begin{itemize}
    \item[(a)] Both facets $A$ and $B$ contain rooks in the first $p-1$ columns  $(p\geq 1)$
    in the priority order $\succ$ (described in (\ref{E:order})) precisely at the same squares
    (positions). Moreover, the facet $B$ contains a rook in the $p$-th column with
    respect to this order and $A$ does not have a rook in the
    $p$-th column.

    \item[(b)] Both facets $A$ and $B$ contain rooks in the first $p-1$ columns in
    the priority order (\ref{E:order}) at the same squares, both $A$ and $B$ contain
    a rook in the $p$-th column and the rook of $A$ is above the rook of
    $B$.
    In other words, if the $p$-th column in the order (\ref{E:order}) is
    $x$, then $x\in A_i$, $x\in B_j$ for some $1<j<i$.
    \item[(c)] Both facets $A$ and $B$ contain rooks in the first $p-1$
    columns $(p\geq 1)$ in the order (\ref{E:order}) at the same squares;
    neither $A$ nor $B$  contains a rook in the $p$-th column and $A'\ll B'$
    where $A'$ and $B'$ are obtained by removing rooks from $A$ and $B$,
    by the following procedure.

    Let $X\subset [m]\times [n]$ be the union of the first $p$
    columns in the order (\ref{E:order}), where by construction $\{p\}\times [n]\subset X$
    is the `empty column' (for both $A$ and $B$). Define $j_i$ (where $0\leq j_i\leq k_i$) as the
    number of rooks in $X\cap \hat{A}$ in the $i$-th row. Remove from the
    chessboard $[m]\times [n]$:
 \begin{enumerate}
 \item[(i)] the `small chessboard' $X$;
 \item[(ii)] the union $Y$ of the first row and all rows where $j_i = k_i$;
 \item[(iii)] the set $Z = A_1\times [n]$.
 \end{enumerate}

Simplices $A'$ and $B'$ are precisely what is left in $A$ and $B$
respectively after the removal of these rooks. Let $j$ be the
number of rows where the equality $j_i = k_i$ is satisfied, that
is $j$ is the number of removed rows aside from the first row.

Note that canonically $K\cong
\Delta^{k_2-j_2,k_3-j_3,\ldots,k_n-j_n;1}_{m-k_1-p, n-j-1}$ and
that the inequality,
\begin{equation}\label{eqn:uslov}
m-k_1-p \geq (k_2-j_2)+\ldots +(k_n-j_n)+ n-j-2
\end{equation}
is a consequence of the inequality $m\geq k_1+k_2+\cdots+k_n+n-1$
and the relation $p-1 = j_2+\ldots + j_n$.
\end{itemize}

The fact that $\ll$ is a linear order is established similarly as
in the proof of Lemma~\ref{lema:linear}.

{\small
\begin{rem}\label{rem:sinisa} \rm
The relation $A'\ll B'$ (in part (c)) refers to the order $\ll$
among facets of the induced multiple chessboard complex $K\cong
\Delta^{k_2-j_2,k_3-j_3,\ldots,k_n-j_n;1}_{m-k_1-p, n-j-1}$. This
isomorphism arises from the canonical isomorphism of chessboards
$([m]\times [n])\setminus (X\cup Y\cup Z) \cong [m']\times [n']$.
As it will turn out in the continuation of the proof there is some
freedom in choosing this isomorphism. Indeed, for each choice of
$A_1$ and $X$ (corresponding to some (c)-scenario) we are allowed
to use the shelling order $\ll$  arising from an arbitrary
isomorphism $K\cong
\Delta^{k_2-j_2,k_3-j_3,\ldots,k_n-j_n;1}_{m-k_1-p, n-j-1}$.
\end{rem}
}

We continue with the proof that $\ll$ is a shelling order for the
multiple chessboard complex
$\Delta^{k_1,\ldots,k_n;\mathbf{1}}_{m,n}$.\footnote{Some readers
may find it convenient to preliminary analyze the special case of
the complex $\Delta_{m,2}^{2,1;\mathbf{1}}$ outlined in
Example~\ref{exam:2-case}.} We are supposed to show that for each
pair of facets $A\ll B$ of
$\Delta^{k_1,\ldots,k_n;\mathbf{1}}_{m,n}$ there exists a facet
$C\ll B$ and a vertex $v\in B$ such that $$A\cap B\subset C\cap
B=B\setminus\{v\}.$$ This statement is established by induction on
$n$ (the number of rows). The claim is obviously true for $n=1$.
Assume that $\ll$ is a shelling order of facets of
$\Delta^{k_1,\ldots,k_q;1}_{d,q}$ for all $q<n$ and all
$k_1,k_2,\ldots,k_q$ such that $d\geq k_1+\cdots+k_q+q-1$.

\medskip
Given the facets $A=(A_1,A_2,\ldots,A_n)$ and
$B=(B_1,B_2,\ldots,B_n)$ such that $A\ll B$, let us consider the
following cases.

\medskip\noindent
$\mathbf{1^\circ}$ Assume that $A_1=\{a_1,a_2,\ldots,a_{k_1}\}<
B_1=\{b_1,b_2,\ldots,b_{k_1}\}$ (in the anti-lexicographic order)
and let $b_i=max(A_1\bigtriangleup B_1)\in B_1$. Observe that in
this case there exists $a_j\in A_1$, $a_j\notin B_1$, $a_j<b_i$.
If there exists a column indexed by $b_0<b_i$ without a rook from
$B$ we let $B'_1=B_1\setminus\{b_i\}\cup\{b_0\}$. Define
$C=(B'_1,B_2,\ldots,B_n)$ and observe that $C\ll B$ and
$$A\cap B\subset C\cap B=B\setminus\{(b_i,1)\}.$$

If all of the columns $1,2,\ldots,b_i$ contain a rook from $B$,
then $a_j\in B_s$ for some $s\in \{2,3,\ldots,n\}$. Let $b_0$ be
the last column in the order (\ref{E:order}) that does not contain
a rook from $B$ and let $B'_s=B_s\setminus\{a_j\}\cup\{b_0\}$. For
the facet $B'=(B_1,\ldots ,B_{s-1},B'_s,B_{s+1},\ldots ,B_n)$ we
have $B'\ll B$ and
$$A\cap B\subset B'\cap B=B\setminus\{(a_j,s)\}.$$

\medskip\noindent
$\mathbf{2^\circ}$ If $A_1=B_1$ we have one of the following
possibilities:

\begin{itemize}
    \item[(a)] $A$ and $B$ contain rooks in the first $p-1$ columns in
    the order (\ref{E:order}) at the same squares; $B$ contains a rook in
    the $p$-th column in (\ref{E:order}) and $A$ does not have a rook in this column.

     Let $b$ denote the label of the $p$-th column, so we have
    that $b\in B_s$ for some $s>1$ and $b\notin A_i$ for all
    $i=1,2,\ldots,n$. Let $b'$ denote the last column in the
    order (\ref{E:order}) that does not contain a rook form $B$,
    and let $B'_s=B_s\setminus\{b\}\cup\{b'\}$.
    For the facet $B'=(B_1,\ldots ,B_{s-1},B'_s,B_{s+1},\ldots ,B_n)$
we have $B'\ll B$ and
$$A\cap B\subset B'\cap B=B\setminus\{(b,s)\}.$$

    \item[(b)] $A$ and $B$ contain rooks in the first $p-1$ columns in
    the order (\ref{E:order}) at the same squares; both $A$ and $B$ contain
    a rook in the $p$-th column and the rook in $A$ is above the rook
    in $B$.

    Again, let $b$ denote the label of the $p$-th column. We have
    that $b\in B_s$ and $b\in A_t$ for some $1<s<t$.
    Let $b'$ be the label of the last column in the
    order (\ref{E:order}) that does not contain a rook from $B$,
    and let $B'_s=B_s\setminus\{b\}\cup\{b'\}$.
    For the facet $B'=(B_1,\ldots ,B_{s-1},B'_s,B_{s+1},\ldots ,B_n)$
we have $B'\ll B$ and
$$A\cap B\subset B'\cap B=B\setminus\{(b,s)\}.$$
    \item[(c)] Both $A$ and $B$ contain rooks in the first $p-1$ columns in
    the order (\ref{E:order}) at the same squares; neither $A$ nor
    $B$  contains a rook in the $p$-th column and $A'\ll B'$ in $K\cong
\Delta^{k_2-j_2,k_3-j_3,\ldots,k_n-j_n;1}_{m-k_1-p, n-j-1}$. By
the inductive assumption (taking into account the inequality
(\ref{eqn:uslov})) we know that $K$ is shellable. Hence, there
exists a facet $C'\ll B'$ of $K$ and a vertex $v\in B'$ such that
$$A'\cap B'\subset C'\cap B'=B'\setminus \{v\}.$$ If we return the
rooks, previously removed from $A$ and $B$ (from the removed rows
and columns of the original complex) and add them to $C'$, we
obtain a facet $C$ of $\Delta^{k_1,\ldots,k_n;\mathbf{1}}_{m,n}$
such that $C\ll B$ and,
$$A\cap B\subset C\cap B=B\setminus \{v\}.$$
\end{itemize}
This observation completes the proof of
Theorem~\ref{thm:shelling-main}. \hfill $\square$

\begin{exam}\label{exam:2-case}{\rm
Here we review the definition of the linear order $\ll$
(introduced in the proof of Theorem~\ref{thm:shelling-main}) and
briefly outline the proof that it is a shelling order for the
complex $\Delta_m :=\Delta_{m,2}^{1,2; \mathbf{1}}$ if $m\geq 4$.

\medskip
By definition the facets of $\Delta_m$ can be described as pairs
$\alpha = (a_1, A_2)$ where $a_1\in [m]$ and $A_2$ is a
$2$-element subset of $[m]$ such that $a_1\notin A_2$. More
explicitly the associated facet is $\hat\alpha = \{(a_1,1)\}\cup
(A_2\times \{2\})\subset [m]\times [2]$.

Let $\mathcal{A}_a = \{(a_1, A_2) \mid a_1 = a\}$  be the
collection of all facets with the prescribed element $a\in [m]$ in
the first row. By definition if $a_1 < b_1$ then $(a_1, A_2)\ll
(b_1, B_2)$ or in other words $\mathcal{A}_{a_1}\ll
\mathcal{A}_{b_1}$ (as sets).

Note that $\mathcal{A}_{a} = \mathcal{A}_{a}^{(1)}\cup
\mathcal{A}_{a}^{(2)} \cup \mathcal{A}_{a}^{(3)}$ (disjoint union)
where,

\begin{equation}\label{eqn:3-slucaja}
\begin{array}{ccl}
\mathcal{A}_{a}^{(1)} & = & \{(a, A) \mid a-1\notin A \}   \\
\mathcal{A}_{a}^{(2)} & = & \{(a, A) \mid a-1\in A, a-2\notin A \} \\
\mathcal{A}_{a}^{(3)} & = & \{(a, A) \mid a-1\in A, a-2\in A \}\\
\end{array}
\end{equation}
(by definition $a-j$ is the unique element $x\in [m]$ such that
$a-j\equiv x \,\mbox{\rm mod}\, (m)$).

By inspection of the definition of $\ll$ we observe that,
\[
\mathcal{A}_{a}^{(1)}\ll \mathcal{A}_{a}^{(2)} \ll
\mathcal{A}_{a}^{(3)}
\]
in the sense that $x_1\ll x_2\ll x_3$ for each choice $x_i\in
\mathcal{A}_{a}^{(i)}$.  It follows from (\ref{eqn:3-slucaja})
that $\mathcal{A}_{a}^{(3)} = \{(a,\{a-1, a-2\})\}$ is a
singleton. The order $\ll$ inside $\mathcal{A}_{a}^{(1)}$ is
determined by reduction to a smaller chessboard (isomorphic to
$[m-2]\times \{1\} \cong [m-2]$). The order $\ll$ inside
$\mathcal{A}_{a}^{(2)}$ is (in the case (c)) determined by
reduction to a smaller chessboard (isomorphic to $[m-3]$).

\medskip
By using this analysis the reader can easily check (case by case)
that $\ll$ is indeed a shelling order on $\Delta_m
:=\Delta_{m,2}^{1,2; \mathbf{1}}$.

\begin{figure}[thb]
\centering
\includegraphics[scale=0.80]{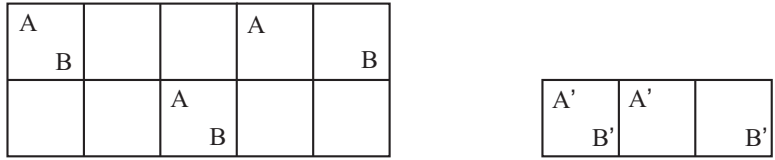}
\caption{Comparison of facets in $\Delta_{5,2}^{1,2;\mathbf{1}}$.}
\label{fig:sah-1}
\end{figure}
\noindent For illustration let $A = (3, \{1,4\})$ and $B = (3,
\{1,5\})$ be two facets of the complex
$\Delta_{5,2}^{1,2;\mathbf{1}}$ (Figure~\ref{fig:sah-1} on the
left). These two simplexes have the same element in the first row.
Moreover the first column in the (lacunary) order (\ref{E:order})
is empty for both $A$ and $B$ so it is the case (c) of the general
comparison procedure that applies here. By removing the first row,
the third column and the second (empty) column we obtain a
$[3]\times [1]$ chessboard and the facets $A'$ and $B'$. This way
the comparison of $A$ and $B$ is reduced to the comparison of $A'$
and $B'$ in $\Delta_{3,1}^{2;\mathbf{1}}$. The most natural choice
is $A'\ll B'$ however (in light of Remark~\ref{rem:sinisa}) we are
free to choose any shelling order on
$\Delta_{3,1}^{2;\mathbf{1}}$.

 }
\end{exam}

\section{An application}
\label{sec:application}

The general colored Tverberg problem, as outlined in
Sections~\ref{sec:intro-1} and \ref{sec:general-color}, is the
question of describing conditions which guarantee the existence of
large intersecting families of multicolored (rainbow) simplices.
By definition a simplex is multicolored if no two vertices are
colored by the same color.

\medskip
We can modify the problem by allowing multicolored simplices to
contain not more than $p$ points of each color where $p\geq 1$ is
prescribed in advance. Following the usual scheme
(Section~\ref{sec:general-scheme}) we arrive at the generalized
chessboard complexes $\Delta_{t,r}^{p,1}$.

\medskip
This connection was one of the reasons why we were interested in
the connectivity properties of multiple chessboard complexes
$\Delta_{t,r}^{p,1}$ and the following statements illustrate some
of possible applications.\footnote{It was kindly pointed by
G\"{u}nter Ziegler that Theorem~\ref{thm:application} is implicit
in \cite{BFZ}, see their Theorem~4.4 and the remark following the
proof of \cite[Lemma~4.2.]{BFZ}.} Further development of these
ideas and new applications to Tverberg-van Kampen-Flores type
theorems can be found in \cite{jvz}.

\begin{theo}\label{thm:application}
Let $r=p^\alpha$ be a prime power. Given $k$ finite sets of points
in $\mathbb{R}^d$ (called colors), of $(p+1)r-1$ points each, so
that $prk\geq (r-1)(d+1)+1$, it is possible to divide the points
in $r$ groups with at most $p$ points of the same color in each
group so that their convex hulls intersect.
\end{theo}

\medskip\noindent
{\bf Proof:} The multicolored simplices are encoded as the
simplices of the simplicial complex $([\ast ]_{\Delta
(p+1)}^{*((p+1)r-1)})^{*k}$. Indeed these are precisely the
simplices which are allowed to have at most $p$ vertices in each
of $k$ different colors. The configuration space of all $r$-tuples
of disjoint multicolored simplices is the simplicial complex,

$$K=(([\ast ]_{\Delta (p+1)}^{*((p+1)r-1)})^{*k})_{\Delta
(2)}^{*r}.$$

Since the join and deleted join commute, this complex is
isomorphic to,

$$K=(([\ast ]_{\Delta (p+1)}^{*((p+1)r-1)})_{\Delta
(2)}^{*r})^{*k}.$$

If we suppose, contrary to the statement of the theorem, that the
intersection of images of any $r$ disjoint multicolored simplices
is empty, the associated mapping $F : K\rightarrow
(\mathbb{R}^d)^{*r}$ would miss the diagonal $\Delta$. By
composing this map with the orthogonal projection to
$\Delta^{\perp}$, and after the radial projection to the unit
sphere in $\Delta^{\perp}$, we obtain a
$(\mathbb{Z}/p)^\alpha$-equivariant mapping,

$$\tilde F : K\rightarrow S^{(r-1)(d+1)-1}.$$

By Corollary \ref{druga'}, the complex $([\ast ]_{\Delta
(p+1)}^{*((p+1)r-1)})_{\Delta (2)}^{*r}$ is $(pr-2)$-connected,
hence the complex $K$ is $(prk-2)$-connected. By our assumption
$prk-2\geq (r-1)(d+1)-1$, so in light of Dold's theorem \cite{M}
such a mapping $\tilde F$ does not exist. \hfill $\square$
\bigskip

Specializing to the case $pk=d+1$, it is easy to see that we could
take $(p+1)r-1$ points of each of $k$ colors in
$\mathbb{R}^{pk-1}$ and obtain the following.

\begin{cor} Let $r$ be a prime power.
Given $k$ finite sets of points in $\mathbb{R}^{pk-1}$ (called
colors), of $(p+1)r-1$ points each, it is possible to divide the
points in $r$ groups with at most $p$ points of the same color in
each group so that their convex hulls intersect.
\end{cor}

Of course, the continuous (non-linear) versions of the above
results are true as well, and with the same proof.

If $r+1$ is a prime, there is a simple proof of this corollary
using the result of \cite{bmz} which even obtains better estimate
($pr$ instead of $(p+1)r-1$) on the number of points of each
color. Namely, one could divide each color of $pr$ points in $p$
"subcolors" of $r$ points each, and obtain the desired division,
even with some additional requirements on the points of the same
color in each group. However, this argument works only in this
case when $r+1$ is a prime.

\begin{rem}\label{rem:generality}{\rm
A more general result related to Theorem~\ref{thm:application} can
be formulated if we allow each of $r$ multicolored sets to contain
$p_1$ points of the first color, $p_2$ points of the second color,
etc. $p_k$ points of the $k$-th color. In this case we arrive at
the complex $\Delta_{t,r}^{p_1,...,p_k;{\mathbf 1}}$ and its
connectivity properties established by Theorem~\ref{thm:main} can
be used again. We omit the details since our main goal in this
paper was to establish improved bounds on the connectivity of
generalized chessboard complexes and Theorem~\ref{thm:application}
was useful to illustrate their importance and to show how they
naturally appear in different mathematical contexts.}
\end{rem}


\begin{thebibliography}{1000}

\bibitem[A04]{A04} C.~Athanasiadis, Decompositions and connectivity of matching
and chessboard complexes, \textit{Discrete Comput. Geom.} 31
(2004), 395–-403.

\bibitem[BL]{BL} I.~B\'{a}r\'{a}ny, D.G.~Larman. A colored version
of Tverberg's theorem. \textit{J.~London Math.\ Soc.}, II.\ Ser.,
45:314--320, 1992.

\bibitem[BSS]{bss} I.~B\'{a}r\'{a}ny, S.B.~Shlosman, and A. Sz\" ucs. On a
topological generalization of a theorem of Tverberg. \textit{J.
London Math. Soc.}, 23:158--164, 1981.

\bibitem[Bj{\"o}95]{bj}
A.~Bj{\"o}rner.
\newblock Topological methods.
\newblock In R.~Graham, M.~Gr\"{o}tschel, and L.~Lov\'{a}sz, editors, {\em
  Handbook of combinatorics}, pages 1819--1872. North Holland, Amsterdam, 1995.

\bibitem[BLV{\v{Z}}]{blvz}
A.~Bj{\"o}rner, L.~Lov{\'a}sz, S.T. Vre{\'c}ica, and R.T.
{\v{Z}}ivaljevi{\'c}.
\newblock Chessboard complexes and matching complexes.
\newblock {\em J. London Math. Soc. (2)}, 49(1):25--39, 1994.

\bibitem[BFZ]{BFZ} P.V.M. Blagojevi\'c, F. Frick, G.M. Ziegler,
{Tverberg plus constraints}, {\em Bulletin of the London
Mathematical Society},  2014, Vol.~46, {953--967}.

\bibitem[BMZ]{bmz}
P.V.M.~Blagojevi{\'c}, B.~Matschke, G.M.~Ziegler.
\newblock Optimal bounds for the colored Tverberg problem.
\newblock {\em J. European Math. Soc.}, Vol.~17, Issue 4, 2015, pp. 739–-754.

\bibitem[BMZ-2]{BMZ-2} P.V.M.~Blagojevi\' c, B.~Matschke, G.M.~Ziegler.
Optimal bounds for a colorful Tverberg--Vre\' cica problem, {\em
Advances in Math.},  Vol. 226, 2011, 5198--5215;
arXiv:0911.2692v2 [math.AT].

\bibitem[FH98]{FH} J. Friedman and P. Hanlon. On the Betti numbers of chessboard
complexes. \textit{J. Algebraic Combin.}, 8:193–203, 1998.

\bibitem[G79]{Garst}  P.~F.~Garst. Cohen-Macaulay complexes and group actions. PhD
thesis, University of Wisconsin-Madison, 1979.

\bibitem[JVZ-2]{jvz} D. Joji\' c, S.T. Vre\' cica, R.T. \v Zivaljevi\' c.
Symmetric multiple chessboard complexes and a new theorem of
Tverberg type,  arXiv:1502.05290 [math.CO].

\bibitem[J08]{Jo-book} J.~Jonsson. \textit{Simplicial Complexes of
Graphs}. Lecture Notes in Mathematics, Vol.\ 1928, Springer 2008.

\bibitem[KRW]{krw}
D.B. Karaguezian, V. Reiner, M.L. Wachs.
\newblock Matching Complexes, Bounded Degree Graph Complexes,
and Weight Spaces of GL -Complexes.
\newblock {\em Journal of Algebra} 239:77--92, 2001.

\bibitem[Koz]{Koz} D.~Kozlov. \textit{Combinatorial Algebraic
Topology}. Springer 2008.

\bibitem[M03]{M} J.~Matou\v sek. \textit{Using the Borsuk-Ulam Theorem.
Lectures on Topological Methods in Combinatorics and Geometry}.
Universitext, Springer-Verlag, Heidelberg, 2003.

\bibitem[SW07]{ShWa} J. Shareshian and M. L. Wachs. Torsion in the matching complex
and chessboard complex. \textit{Adv. Math.}, 212(2):525–-570,
2007.

\bibitem[T66]{Tve} H.~Tverberg. \newblock A generalization of Radon's theorem. \textit{\em J.
London Math.\ Soc.}, 41:123--128, 1966.

\bibitem[V\v Z94]{VZ94} S.~Vre\' cica and R.~\v Zivaljevi\' c. New cases of
the colored Tverberg theorem. In H.~Barcelo and G.~Kalai, editors,
\textit{Jerusalem Combinatorics '93}, Contemporary Mathematics
Vol.\ 178, pp. 325--334, A.M.S.\ 1994.

\bibitem[V{\v{Z}}11]{vz11}
S.~Vre{\'c}ica, R.~{\v{Z}}ivaljevi{\'c}.
\newblock Chessboard complexes indomitable.
\newblock {\em J. Combin. Theory Ser. A}, 118(7):2157--2166, 2011.

\bibitem[Zi94]{zie} G.M.~Ziegler. \newblock Shellability of chessboard
complexes. {\em Israel J. Math.} 1994, Vol.~87, 97--110.


\bibitem[Zi11]{Ziegler}  G.M.~Ziegler.  \newblock $3N$ colored points in a plane.
\newblock {\em Notices of the A.M.S.} Vol.\ 58 , Number 4,
550--557, 2011.

\bibitem[\v Z96]{User1}
R.~\v Zivaljevi\' c.
\newblock User's guide to equivariant methods
in combinatorics.  \textit{Pu\-bli\-cations de l'Institut
Mathematique} (Beograd), 59(73), 114--130, 1996.

\bibitem[\v Z98]{User2}
R.~\v Zivaljevi\' c.
\newblock User's guide to equivariant methods
in combinatorics II. \textit{Publi\-cations de l'Institut
Mathematique} (Beograd), 64(78) 1998, 107--132.

\bibitem[\v Z99]{TV-bundle} R.~\v Zivaljevi\'{c}.
The Tverberg-Vre\'{c}ica problem and the combinatorial geometry on
vector bundles. \textit{Israel J.\ Math}, 111:53–-76, 1999.

\bibitem[\v Z04]{Z04}
R.T. \v Zivaljevi\'{c}. Topological methods. Chapter 14 in
\textit{Handbook of Discrete and Computational Geometry}, J.E.\
Goodman, J.\ O'Rourke, eds, Chapman \& Hall/CRC 2004, 305--330.


\bibitem[{\v{Z}}V92]{zv92}
R.T. {\v{Z}}ivaljevi{\'c} and S.T. Vre{\'c}ica.
\newblock The colored {T}verberg's problem and complexes of injective
  functions.
\newblock {\em J. Combin. Theory Ser. A}, 61(2):309--318, 1992.


\end{thebibliography}
\end{document}